# 九年级学生数学运算能力培养的实践研究

A Practical Study on Developing Mathematical Computation Ability of Ninth-Grade Students

刘杨，北京大学附属中学

LIU Yang, The Affiliated High School of Peking University

李雪丹，北京市顺义牛栏山第一中学

LI Xuedan, The Beijing Shunyi Niulanshan No.1 Middle School

内容提要：本研究结合《义务教育数学课程标准（2011 年版）》及《普通高中数学课程标准（2017 年版 2020 年修订）》，参考了 PISA 测评理论的部分内容，确定了本次实践研究着力通过算理分析和算法对比两个方向提升九年级学生的数学运算素养水平。选取笔者所在学校九年级的学生共 171 人参与实践研究，在 2021 年 9 月至 12 月期间参与了持续的针对数学运算的讨论课程。2022 年 1 月使用全区统一调研试卷作为学生数学运算能力素养水平的评价工具进行了相关测试。收集到有效数据 171 份，之后借助 Excel 365 及 SPSS 23.0 对测评得到的数据进行了整体分析、相关性分析、测评维度分析。

根据实践及测评结果，本文提出以下三个建议帮助教师提升学生的数学运算素养水平：（1）激发学生学习运算的兴趣；（2）帮助学生更好地明晰算理；（3）提升学生的算法选择能力。

主题词：数学运算　能力培养　实践研究　完成时间：2022 年 3 月

# Abstract

This study conducted practical research on cultivating students' mathematical operation ability through literature analysis, theoretical practice, and collection and statistical analysis of data. The author first collected and analyzed relevant research background literature to determine the development status of related practical research. Based on the "Mathematics Curriculum Standards for Compulsory Education (2011 Edition)" and the "Mathematics Curriculum Standards for Ordinary High Schools (2017 Edition Revised in 2020)", and referring to some content of the PISA assessment theory, the study focused on improving the mathematical operation literacy level of ninth-grade students through two directions: algorithmic analysis and algorithm comparison.

A total of 171 ninth-grade students from the author's school participated in the practical research, with Experimental Group A (57 students), Experimental Group B (58 students), and Control Group C (56 students). The students in the two experimental groups participated in continuous discussion courses on mathematics operations from September to December 2021. In January 2022, the unified research paper was used as the evaluation tool for students' mathematical operation literacy level during the final district-wide research. Valid data was collected from 171 students, and the evaluation data was analyzed comprehensively, correlatively, and dimensionally using Excel 365 and SPSS 23.0.

The results of the evaluation showed that: (1) The initial mathematical operation ability level of the experimental and control group students was basically the same, with no significant differences, but after a long period of practice, the students' operation ability level improved. (2) The degree to which students' mathematical operation ability level was improved through continuous targeted discussion and research was not significantly related to gender but had a moderate positive correlation with the number of times the students participated in related discussions. (3) Most ninth-grade students can reach the level one requirements of mathematical operation literacy in high school curriculum standards, and continuous analysis and comparison of algorithms and logic can to some extent improve students' ability to achieve level two requirements. Further research is needed to reach level three requirements.

Based on the practical research and evaluation results, this study proposes several suggestions to help teachers improve students' mathematical operation literacy level: (1) Stimulate students' interest in learning operations, mobilize their subjective initiative, and encourage their participation by moderately stimulating students' exploration enthusiasm with typical problems. (2) The improvement of mathematical operation ability cannot rely solely on practicing problems but also needs to focus on selecting more typical problems to help students clarify algorithms. (3) While clarifying algorithms, attention should also be paid to methods to cultivate students' good mathematical thinking qualities, and improve their algorithm selection ability by analyzing problems and choosing optimal operation methods.

# 九年级学生数学运算能力培养的实践研究

**提要**　本研究主要通过文献分析，理论实践、收集统计分析等方法对学生数学运算能力的培养进行了实践研究。笔者首先通过对相关研究背景的文献收集与分析确定了相关实践研究的发展现状，结合《义务教育数学课程标准（2011年版）》及《普通高中数学课程标准（2017年版2020年修订）》，参考了PISA测评理论的部分内容，确定了本次实践研究着力通过算理分析和算法对比两个方向提升九年级学生的数学运算素养水平。选取笔者所在学校九年级的学生共171人参与实践研究，设立实验组A（57人）、实验组B（58人）和对照组C（56人），两个实验组的学生在2021年9月至12月期间参与了持续的针对数学运算的讨论课程。2022年1月借助期末全区统一调研的机会，使用统一调研试卷作为学生数学运算能力素养水平的评价工具进行了相关测试。收集到有效数据171份，之后借助Excel 365及SPSS 23.0对测评得到的数据进行了整体分析、相关性分析、测评维度分析。

相关测评结果表明：（1）实验组和对照组学生数学运算能力水平初始状态基本一致，无显著性差异，通过一段较长时间的实践，可以看到学生的运算能力水平却有提升；（2）通过持续的针对性讨论研究提升学生的数学运算能力水平的程度与学生的性别没有显著相关，与学生参与相关讨论的次数有中等程度的正相关性；（3）九年级学生经过学习大部分人可以达到高中课程标准中数学运算素养水平一的要求，而通过持续的分析算理和算法对比学习可以在一定程度上提升学生数学运算素养水平二的能力，对于水平三的要求如何达成还需要进一步研究。

根据实践及测评结果，本文提出以下几个建议帮助教师提升学生的数学运算素养水平：（1）激发学生学习运算的兴趣，调动学生的主观能动性，激发学生的参与意识，以典型问题适度激发学生的探讨热情，都有助于让学生主动参与数学运算的学习，提高学生学习运算的兴趣；（2）数学运算能力的提升不能仅依靠练习题目，还需要在选题上多下功夫，让题目更具有典型性，可以帮助学生更好地明晰算理；（3）在明晰算理的同时要讲究方法着力培养学生良好的数学思维品质，通过分析题目选择更优的运算思路，设计合理的运算程序，提升学生的算法选择能力。

## 一、研究背景
### （一）　数学运算素养培养的目标要求贯穿基础教育

在《义务教育数学课程标准（2011年版）》中，数学运算能力就已经被列为描述学生数学能力的十个核心词之一，之后的《普通高中数学课程标准（2017年版2020年修订）》中又将数学运算列为了六大数学核心素养之一，这表明对学生运算素养的培养、提升和评价是贯穿在整个基础教育过程中的。

在《义务教育数学课程标准（2011年版）》中认为"运算能力主要是指能够根据法则和运算律正确地进行运算的能力。培养运算能力有助于学生理解运算的算理，寻求合理简洁的运算途径解决问题。"[1]

在《普通高中数学课程标准（2017年版2020年修订）》中将数学运算列为六大数学核心素养之一，其中将数学运算这一核心素养界定为"在明晰运算对象的基础上，依据运算法则解决数学问题的素养。"[2]

显然在这一层面上，义务教育阶段的提出的"运算能力"和高中阶段课标中界定的"数学运算素养"在核心内涵上是一脉相承的。

当然高中阶段明确指出学生要"通过高中数学课程的学习，进一步发展数学运算能力；有效借助运算方法解决实际问题；通过运算促进数学思维发展,形成规范化思考问题的品质，

养成一丝不苟、严谨求实的科学精神。"在这个层面上，处于义务教育阶段的学生显然需要更多的发展和成长。

其中九年级作为义务教育阶段的最后一年，更应该着力于发展学生的数学运算能力，向高中阶段要求的数学运算素养标准靠拢，从而为学生的成长奠定更坚实的基础。

### （二） 运算素养培养的实践研究有待深入

通过中国知网论文发布数据趋势分析功能可以发现（图1-1），自《普通高中数学课程标准（2017年版）》发布以来，国内以运算素养为主题的相关研究论文从2018年的179篇迅速增长到2019年的330篇，并持续保持相关研究热度，2021年发表相关论文320篇。由此可见对运算素养的研究在国内基础教育领域受到了持续关注。

图1-1

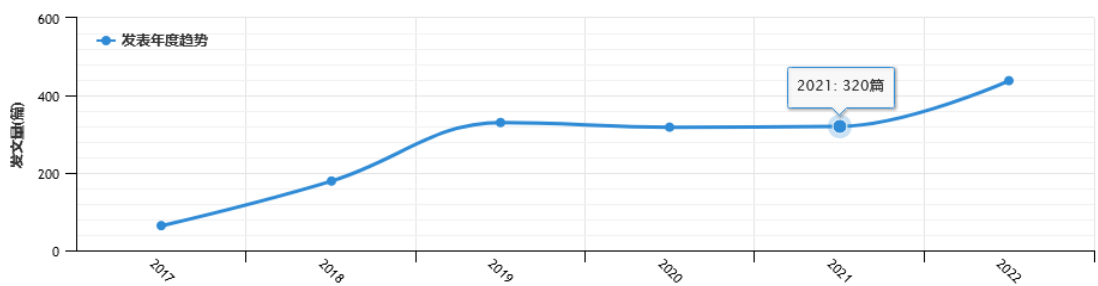

深入分析2021年的320篇论文可以发现其中除学位论文70篇之外以"核心素养"作为主要主题筛选时有129篇论文，而分别以"数学运算"、"运算素养"、"运算能力"、"数学运算素养"为主要主题进行筛选，只能在2021年检索到10~80篇不等的结果，这表明虽然有大量的论文可以用"运算素养"作为主题被检索到，但在论文中并没有将"运算素养"作为主要的研究主题，而是在较为宽泛地研究"核心素养"这个整体，真正聚焦到运算素养的研究论文相对较少。

当进一步以"运算素养"、"培养"、"实践研究"作为主题检索时，在2021年共检索到相关论文8篇，其中5篇为硕士毕业论文，均为研究高中运算素养的相关研究，其余论文为发表于各级期刊由一线教师开展的相关研究，涉及高中学段两篇、小学一篇，对于初中阶段运算素养的培养未有研究涉及。

图1-2

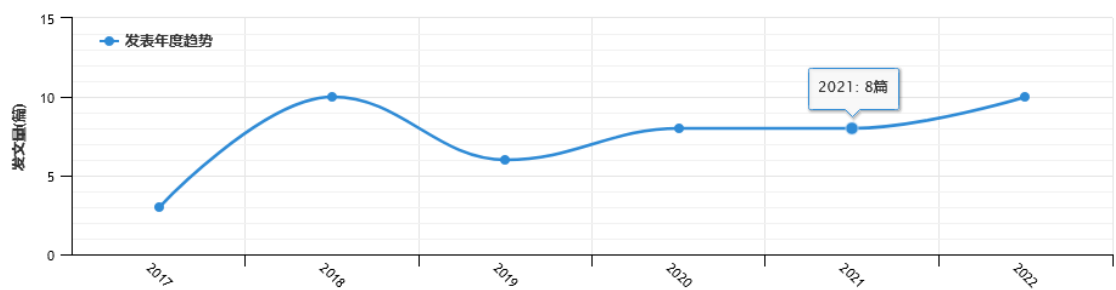

由统计图1-2可见，对于数学运算素养培养的实践研究虽然呈现上升趋势，但相关的研究更多由高校研究生主持开展，一线教育工作者参与较少。且相关研究的持续时间往往受限于研究者的学业，很难长时间开展，相关研究有待进一步深入。

### （三） 一线教师开展持续性的实践研究具有一定优势

很多的实践研究需要在较长的时间段上展开持续的追踪观察，基于这一实践研究的特点，

处于教学一线教师更适合开展相关工作，并在实践过程中得到第一手的数据。此类研究的难点主要在于构建相关的理论去设计相应的实践活动，并指导实践研究的开展。

由于对运算素养的相关理论研究已经开展的比较充分，对于如何培养学生的数学运算素养已经产生了一些初步的成果，故一线教师可以基于对此类研究的学习，构建相关的实践活动框架，并完成相关的实践研究的工作。

本文笔者就是基于此开展的相关实践研究工作。

## 二、相关理论
### （一）相关界定
1. 《普通高中数学课程标准（2017 年版 2020 年修订)》中对数学运算素养水平的界定

在《普通高中数学课程标准（2017 年版 2020 年修订)》中对数学运算素养水平给出了三个水平的界定，具体要求如表 2-1 所示。

表 2-1

| 水平 | 数学运算 |
| --- | --- |
| 水平一 | 能够在熟悉的数学情境中了解运算对象，提出运算问题。<br>能够了解运算法则及其适用范围，正确进行运算；能够在熟悉的数学情境中，根据问题的特征形成合适的运算思路，解决问题。<br>在运算过程中，能够体会运算法则的意义和作用，能够运用运算验证简单的数学结论。<br>在交流的过程中，能够用运算的结果说明问题。 |
| 水平二 | 能够在关联的情境中确定运算对象，提出运算问题。<br>能够针对运算问题，合理选择运算方法、设计运算程序，解决问题。<br>能够理解运算是一种演绎推理；能够在综合运用运算方法解决问题的过程中，体会程序思想的意义和作用。<br>在交流的过程中，能够借助运算探讨问题。 |
| 水平三 | 在综合的情境中，能够把问题转化为运算问题，确定运算对象和运算法则，明确运算方向。<br>能够对运算问题，构造运算程序，解决问题。<br>能够用程序思想理解和表达问题，理解程序思想与计算机解决问题的联系。<br>在交流的过程中，能够用程序思想理解和解释问题。 |

2. 国际学生评估项目（PISA）的相关界定[3]

国际学生评估小项目（PISA）测评体系有三个维度构成，分别是数学内容、数学过程以及数学情境,数学素养测评的理论框架由这三个维度构成的一级指标和其下分的具体分类一道组成，为直观体现上述理论框架特制作下图 2-1 理论框架结构图。

图 2-1 PISA2012 数学素养理论框架结构图

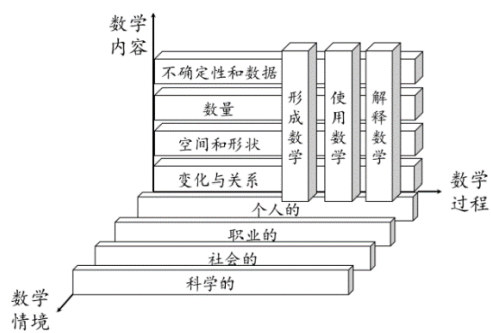

从图 2-1 中可以看到，PISA2012 数学素养理论框架结构总共包含 3 个维度，它们是数学内容维度、数学过程维度以及数学情境维度。当然 PISA 测试的特点在数学情境角度上和我们日常的教学还是有比较大的差别的，能够比较好地和《高中课程标准（2020 年修订）》中运算素养相关内容融合的是其数学内容维度和数学过程维度。

其中数学内容的 4 个方面中"数量"部分与数学运算素养是息息相关的；从数学过程维度来看，其主要涉及的形成数学、使用数学、解释数学 3 个过程可以大致与《高中课程标准（2020 年修订）》中给出的三个水平界定对应。

3. 义务教育阶段提升学生数学运算能力水平的策略

在相关研究著作中 [4,5]，作者提出了多种培养学生数学运算能力的策略，通过对比总结作者认为培养提升学生运算能力，至少可以从以下几个方面着手：

(1) 激发学生学习运算的兴趣；
(2) 增加板书，示范规范的运算流程，培养学生养成良好的运算习惯；
(3) 重视数学基础知识的教学积累，明确算理；
(4) 讲究方法培养学生良好的数学思维品质，合理选择算法；
(5) 增强学生的符号意识；
(6) 适当地合理使用学习辅助工具。

**（二） 本次实践研究的理论框架**

综合前述讨论，结合第一作者所在学校的学情，考虑到九年级学生的学业特征和切实需求，在本次实践研究中，着重从以下几个维度去刻画和描述学生的数学运算素养水平。

1. 从数学运算内容维度描述

九年级学生需要学习的和运算直接相关的内容包括：一元二次方程、二次函数、反比例函数三大部分内容,在学习过程中亦会需要初一年级和初二年级数与代数部分运算内容的支持。故本次实践研究中作为数学运算能力培养教学的知识载体主要是一元二次方程、二次函数和反比例函数三部分内容。

2. 从数学运算过程维度描述

解决数学问题的过程其实是相当复杂的,在解答的过程中的表现体现了学生数学运算素养所能达到的水平。对此《高中课程标准（2020 年修订）》中给出了三个水平描述，但高中课程标准所要求的水平显然是针对高中学生提出的，考虑到九年级学生的发展特征，笔者认为对于这一年龄段的学生来说水平三的要求相对过高，故在本次实践研究中，重点关注对学生达成"水平一"和"水平二"的培养。

3. 数学运算培养策略

在相关实践研究中提出的诸多培养策略，结合九年级的教育教学特点，在本次实践研究中将着重尝试从算理和算法两个角度出发对学生的数学运算素养进行培养。因为算理和算法是数学运算的基础，作为解题基础的算理是发展学生数学运算能力的重要前提。在九学习备考阶段也需要通过复习相关数学概念、运算规律、运算性质等帮助学生再次复习熟悉相关算理，以算理分析作为研究的视角事半功倍。

另一方面算法是运算的方法，在九年级的教育教学中要帮助学生解决"如何计算"的问题，一定不能回避对算法的梳理，在面对数学问题时选择不同的算法，解决问题时的复杂程度也会有所的差异，通过在算法学习中对题目不同解法的对比总结，可以帮助九年级的学生进一步理解算理和算法的关系，并在解题时选择恰当的视角，使用合理的算法去解决问题。

故在本次实践研究中，将重点关注：
(1) 数学运算算理的分析；
(2) 数学运算算法的对比；
并检验从这两个角度出发对学生进行引导，是否可以有效提高学生的数学运算能力，从而提升学生的数学运算的素养水平。

### 三、实践设计
#### （一）实践对象的确定
本次实践研究的对象是从第一作者所在学校 2022 届学生中抽取的，以被抽取的这部分同学参与海淀区数学统一调研的成绩作为研究对象。实际上选择 2022 届毕业年级学生参与区统一调研的成绩作为研究对象有以下几个有利因素。

第一，在确定参与实践研究的学生名单时，第一作者已经在本年级任教一年，熟悉年级中各层次水平学生的实际表现；第二，第一作者随本年级学生进入毕业年级教学，在实践研究过程中作为毕业年级的任课教师，第一作者可以较好地把控实践过程，也便于开展访谈等研究活动；第三，根据相关教学安排，初中毕业年级学生在九年级第一学期将学习二次函数等重要的知识内容，在这个过程中有利于通过实践培养学生在数学运算方面的能力；第四，毕业年级学生有更强的动力去改进自身不足，提升个人解题能力，故参与相关实践活动的积极性更有保障；第五，学生运算能力的评价可以借助编制相对合理的区级九年级统一调研完成。

#### （二）实践流程
参与本次实践研究的学生共分为 2 个实验组和 1 个对照组，共有 180 名同学参与，所有同学均衡分布在年级各行政班级中，无显著分布差异，相关实践流程如下：

图 3-1

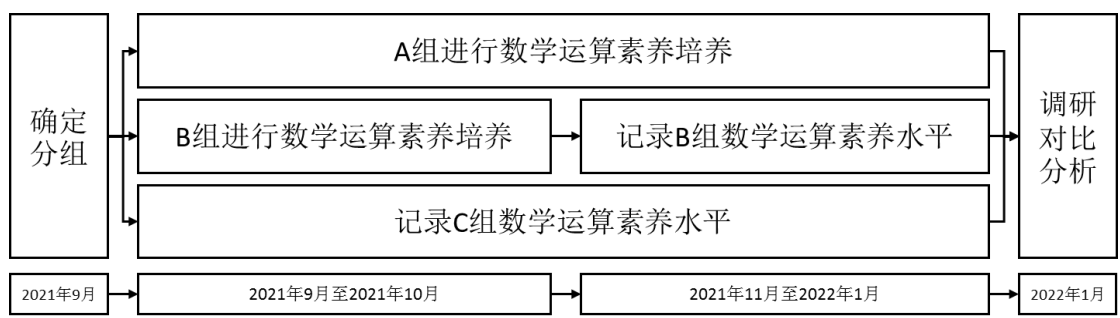

根据图 3-1 所示流程，本次实践研究基本分为四个阶段：
1. 准备阶段：从 2021 年 8 月开始，截止 2021 年 9 月开学第一周，目标为确定实验组和对照组参与学生。确定方法是以所有学生八年级第二学期期中和期末调研成绩标准分的加权平均值（期中占 40%，期末占 60%）按从高到低进行排序，之后采用系统抽样法，筛选出 A、B 两个实验组各 60 位同学以及一个包含 60 名同学的实验对照组 C。之后通过说明会征询相关同学的参与意见，最终确定参与实践研究的学生名单。
2. 实施阶段 1，从 2021 年 9 月中旬开始持续到 2021 年 10 月底，主要实践工作为：
   (1) 实验组 A 和实验组 B 每周各用一天的中午参与实践研究，形式为与运算相关内容的课堂讨论，讨论围绕运算算理分析和运算算法对比展开，包括综合应用的讨论。配合学校九年级教学具体进度选择合适的载体展开讨论；
   (2) 对照组 C 不参与相关实践活动讨论

(3) 收集反映实验组和对照组学生运算能力水平情况的相关调研成绩数据。
3. 实施阶段2：从2021年11月中旬开始持续到2021年12月底，主要实践工作为：
   (1) 实验组A继续参与每周一次的实践研究课，形式仍然为运算相关内容的课堂讨论，继续围绕运算算理分析和运算算法对比展开，包括综合应用的讨论。进度在配合学校九年级教学具体进度的同时，适当选择合适的中考复习内容作为载体展开讨论；
   (2) 实验组B和对照组C不参与这个阶段的实践活动讨论；
   (3) 继续收集反映实验组和对照组学生运算能力水平情况的相关调研成绩数据。
4. 反馈分析阶段：2022年1月，借助2022年学区九年级期末数学调研试卷中直接考察数学运算素养的题目，收集参与实践研究的实验组A、B以及对照组C相关学生的答题数据，对各组学生的数学运算能力素养水平的差异变化进行评价。

### 四、实践与调研过程
#### （一）确定实验组与对照组

2021年9月开学后，根据第一作者所在学校毕业年级学生分班名单，调取了所有学生八年级第二学期期中和期末调研成绩，转换为标准分后进行加权平均计算（期中占40%，期末占60%)得到入围标准，将所有学生按入围标准从高到低进行排序，之后采用系统抽样法，筛选出实验组A、实验组B、对照组C各60名同学作为实践候选参与者，一共180人。所有180名实践参与者皆来自于平行班级，涉及毕业年级所有16个班级，其中男生共101人，女生共79人，相对占比如下图4-1所示。

图4-1

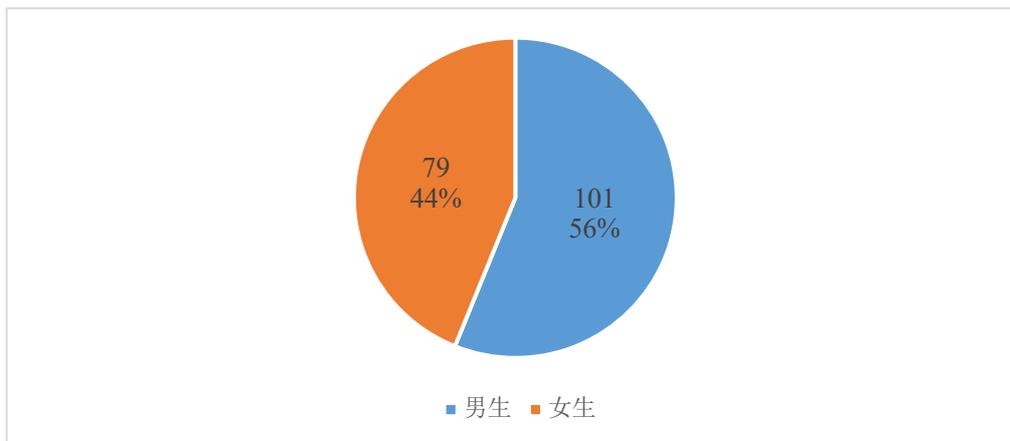

考虑到年级学生总数的男女生占比情况为男生占该年级学生总数的61%，女生占39%，可以认为抽样结果与学生实际比例基本一致。对三个分组中的学生的入围标准分分布情况进行Pearson相关性分析，其中实验组A、实验组B和对照组C分别给出分组编号1、2和3，相关分析结果整理如下表4-1所示。

表4-1 各分组入围标准分相关性分析

|  |  | 各分组 | 入围标准分 |
|---|---|---|---|
| 各分组 | Pearson 相关性 | 1 | -0.100 |
|  | 显著性（双侧） |  | 0.183 |
|  | N | 180 | 180 |
| 入围标准分 | Pearson 相关性 | -0.100 | 1 |
|  | 显著性（双侧） | 0.183 |  |
|  | N | 180 | 180 |

从表 4-1 可知，各分组学生的入围标准分的 Pearson 相关系数为-0.100，显著性（双侧）水平为 0.183。这表明两个实验组和一个对照组学生在入围本次实践研究时，其数学运算素养水平无显著性差异，可以在之后进行观察比较。之后和实验组 A 与实验组 B 的学生进行了沟通确认了相关学生的个人参与意愿，并和学生明确每次相关实践活动都自愿参与，其中实验组 A 候选 60 人中有 8 名同学选择不参与，其中 5 人同意收集不带其个人敏感信息的相关调研数据进行研究，实验组 B 候选 60 人中有 10 名同学选择不参与，其中 8 名同学同意收集不带其个人敏感信息的相关调研数据进行研究，对照组 C 的 60 名同学中有 4 人拒绝提供相关调研数据进行研究，其他 56 名同学授权收集不带个人敏感信息的相关调研数据进行研究。最终参与实践研究的学生情况统计如下表 4-2 所示。

表 4-2 各分组最终参与本研究学生情况统计

| 序号 | 组别 | 参与人数 | 参与率 | 同意收集数据 | 收集率 | 退出人数 |
| --- | --- | --- | --- | --- | --- | --- |
| 1 | 实验组 A | 52 | 87% | 57 | 95% | 3 |
| 2 | 实验组 B | 50 | 83% | 58 | 96% | 2 |
| 3 | 实验组 C | 56 | 93% | 56 | 93% | 4 |

（二）实践过程

本次实践研究的每次讨论课都只允许实验组的学生参与，但并不强制要求实验组的所有学生参与，学生每次的讨论都可以视个人情况自愿参加，只在每次讨论课时强制要求所有到场学生完成签到，各次讨论课的签到记录情况见表 4-3。

表 4-3 各次讨论课个人签到情况汇总统计

| 课次 | 知识板块 | 讨论重点 | 实验组 A 签到人数 | 实验组 B 签到人数 |
| --- | --- | --- | --- | --- |
| 1 | 方程与不等式 | 运算算理分析 | 47 | 41 |
| 2 | 方程与不等式 | 算法思路对比 | 36 | 33 |
| 3 | 方程与不等式综合 | 运算综合 | 32 | 34 |
| 4 | 函数 | 运算算理分析 | 32 | 29 |
| 5 | 函数 | 算法思路对比 | 31 | 27 |
| 6* | 函数综合 | 运算综合 | 33 | 30 |
| 7 | 函数综合 | 运算综合 | 32 | 33 |
| 8 | 数与式 | 运算算理分析 | 36 | — |
| 9 | 数与式 | 算法思路对比 | 33 | — |
| 10 | 方程与不等式 | 算法思路对比 | 33 | — |
| 11 | 方程与不等式综合 | 运算综合 | 35 | — |
| 12 | 函数综合 | 运算综合 | 33 | — |
| 13 | 函数综合 | 运算综合 | 37 | — |

\* 代表附录 1 中具体内容示例对应的课次

由表 4-3 可知，13 次专题讨论内容设计数与式、方程与不等式、函数三个初中阶段涉及学生数学运算素养的主要知识板块展开，根据本年级教学进度配合选取各次讨论课的主要内容。从表中可以看出在本次实践研究的实施阶段 1 中，由于学生可以自由决定是否参与某一次的讨论课，讨论课的参与人数在实验组 A 和实验组 B 中均呈现出不同程度的下降，最后两次讨论（课次 6 和课次 7）的实验组 B 参与人数小幅回升与几名学生交流后确认为 11 月期中考试备考导致，同时在实验组 A 的学生中可以看到在 11 月期中考试之后参与讨论课的人数有更明显的回升，对其中部分同学访谈后确认为感觉讨论课对成绩有帮助，自愿多参与

讨论课活动。13 次专题讨论的具体实践情况如下：

2021 年 9 月开学第一周实验组 A 和实验组 B 开始参与相关实践研究，实验组 B 于 9 月 9 日周四中午 12:30 至 13:10 进行了第一次专题讨论，实验组 A 于 9 月 10 日周五中午 12:30 至 13:10 进行了第一次讨论，每次讨论课两个实验组讨论的主题相同，实施教师为同一教师。根据实践计划，各次专题讨论重点均围绕运算算理分析和运算思路对比展开，直到 2021 年 10 月 29 日周五为止，在九年级期中考试前共进行了 7 次专题讨论。

2021 年 11 月期中考试之后，实验组 B 的学生不再参与后续实践研究，实验组 A 的学生于 11 月 12 日周五中午 12:30 至 13:10 开始继续每周一次的专题讨论，至 2021 年 12 月 24 日共开展专题讨论 6 次（其中 2021 年 12 月 17 日的专题讨论因学生备考 18 日中考英语听说口语第一次测试而取消），之后实验组 A 的运算专题讨论实践停止。

### 五、调研对实践的反馈与分析

2022 年 1 月 17 日获取到学区九年级期末数学调研试卷学生小题得分统计表，根据学区发布的此次期末调研双向细目表及学生实际答题情况，将表 5-1 所列题号对应的问题作为本次实践研究中用于评价学生数学运算素养评价的题目，具体问题见附录 2。

表 5-1 九年级第一学期期末数学调研考察运算能力题目统计

| 题号 | 分值 | 题型 | 试题类型 | 知识板块 | 能力板块 | 素养水平 |
| --- | --- | --- | --- | --- | --- | --- |
| 1 | 2 | 选择题 | 基础性试题 | 函数 | 运算能力 | 数学运算水平一 |
| 13 | 2 | 填空题 | 基础性试题 | 方程与不等式 | 运算能力 | 数学运算水平一 |
| 17 | 5 | 解答题 | 基础性试题 | 方程与不等式 | 运算能力 | 数学运算水平一 |
| 18 | 5 | 解答题 | 基础性试题 | 方程与不等式 | 运算能力 | 数学运算水平一 |
| 19-1 | 3 | 解答题 | 基础性试题 | 函数 | 运算能力 | 数学运算水平一 |
| 19-2[*] | 2 | 解答题 | 基础性试题 | 函数 | 运算能力 | 数学运算水平二 |
| 22-1 | 2 | 解答题 | 基础性试题 | 方程与不等式 | 运算能力 | 数学运算水平一 |
| 22-2[*] | 4 | 解答题 | 基础性试题 | 方程与不等式 | 运算能力 | 数学运算水平二 |
| 26-1 | 1 | 解答题 | 基础性试题 | 函数 | 运算能力 | 数学运算水平一 |
| 26-2[**] | 3 | 解答题 | 发展性试题 | 函数 | 推理能力 运算能力 | 数学运算水平二 |
| 26-3[**] | 2 | 解答题 | 发展性试题 | 函数 | 抽象概括能力 运算能力 | 数学运算水平三 |

\* 标记的问题为根据学生答题情况与题目描述二次划分的结果，区数据中未分开描述
\*\* 标记的题目选择不同的思路视角，可以借助推理和抽象概括规避大量的运算过程

#### （一）整体分析

本次调研直接涉及数学运算能力的题目 7 题 9 问，满分 31 分，整个调研试卷满分 100 分。阅卷时全区统一评分标准，各题设立评分组，采取电脑随机双评机制，阅卷误差为 0 分，当双评出现不同分数时由评分组组长进行仲裁，从而确保了学生分数的准确性。从学区反馈的学生小题分数据中筛选出所有参与实践研究的学生数据，进行整理后可以得到基本信息的统计表 5-2 和表 5-3。

表 5-2 参与实践研究学生本次调研总成绩情况

|  | 均值 | 中位数 | 众数 | 标准差 | 最高分 | 最低分 |
| --- | --- | --- | --- | --- | --- | --- |
| 实验组 A | 82.56 | 84 | 90 | 9.46 | 97 | 43 |
| 实验组 B | 70.91 | 74.75 | 89 | 17.97 | 96.5 | 10 |
| 对照组 C | 73.49 | 80.50 | 83,81 | 18.04 | 99 | 25 |
| 整体情况 | 75.64 | 81 | 89 | 16.38 | 99 | 10 |

表 5-3 参与实践研究学生本次调研运算能力考题得分情况

| | 均值 | 中位数 | 众数 | 标准差 | 最高分 | 最低分 |
|---|---|---|---|---|---|---|
| 实验组 A | 25.51 | 26 | 26 | 3.16 | 31 | 14 |
| 实验组 B | 21.25 | 23 | 25 | 6.57 | 31 | 1 |
| 对照组 C | 22.30 | 25.25 | 26 | 6.82 | 31 | 2 |
| 整体情况 | 23.01 | 25 | 26 | 6.01 | 31 | 1 |

虽然在表 5-2 中从各项基本数据上看在本次调研中实验组 B 的成绩相比对照组 C 的成绩都有明显差距，但通过相关性检验却发现两组之间考试成绩与分组之间并无显著性差异，如表 5-4 所示。

表 5-4 实验组 B 与对照组 C 成绩相关性检验

| | | 分组 | 调研标准分分布 |
|---|---|---|---|
| 分组 | Pearson 相关性 | 1 | 0.084 |
| | 显著性（双侧） | | 0.377 |
| | N | 114 | 114 |
| 调研标准分分布 | Pearson 相关性 | 0.084 | 1 |
| | 显著性（双侧） | 0.377 | |
| | N | 114 | 114 |

实际上通过分析表 5-3 可以进一步看出两组学生在成绩裸分上的差距至少并非是学生运算能力水平的差距导致的，实验组 B 和对照组 C 的学生在运算能力水平上处于同一水平。在此基础上，将所有学生本次测验的成绩均转为标准分，之后展开进一步分析。

（二）相关性分析

1. 分组相关性分析

将实验组 A 纳入分析序列后，对实验组 A（序号 1）、实验组 B（序号 2）和对照组 C（序号 3）进行成绩标准分的 Pearson 相关性分析，得到表 5-5。

表 5-5 各分组调研成绩标准分相关性分析

| | | 各分组 | 调研标准分 |
|---|---|---|---|
| 各分组 | Pearson 相关性 | 1 | -0.205[**] |
| | 显著性（双侧） | | 0.007 |
| | N | 171 | 171 |
| 调研标准分 | Pearson 相关性 | -0.205[**] | 1 |
| | 显著性（双侧） | 0.007 | |
| | N | 171 | 180 |

[**] 在 0.01 水平（双侧）上显著相关

从表 5-5 中可以看出各分组间的成绩与其分组情况在 0.01 水平上显著相关，相关系数为-0.205 属于弱相关。由此可以得出，对教师而言，学生运算能力的培养需要一定的时间持续进行才可能产生一点点效果，这一数学素养的培养并不能毕其功于一役，而是需要持之以恒地进行。对于学生来说这一情况可能更为明显。

2. 参与讨论次数的相关性分析

对于实验组 A 和实验组 B 来说，通过人为限制参与讨论的次数（实验组 A 最多可参加 13 次讨论，实验组 B 最多可参加 7 次讨论）和随机因素确定参与次数（学生自愿决定是否参与某一次讨论），可以得到精确到个人的参与讨论次数的统计情况（表 4-3 为个人参与情况汇总结果），以此为基础分析学生个人参与讨论次数与其期末调研运算能力水平标准分之间的 Pearson 相关性，得到表 5-6A 和表 5-6B。

表 5-6A 实验组 A 学生参与讨论次数与运算能力水平标准分相关性分析

|  |  | 参与次数 | 运算能力标准分 |
|---|---|---|---|
| 参与次数 | Pearson 相关性 | 1 | 0.502** |
|  | 显著性（双侧） |  | 0.000 |
|  | N | 57 | 57 |
| 运算能力标准分 | Pearson 相关性 | 0.502** | 1 |
|  | 显著性（双侧） | 0.000 |  |
|  | N | 57 | 57 |

** 在 0.01 水平（双侧）上显著相关

表 5-6B 实验组 A 学生参与讨论次数与运算能力水平标准分相关性分析

|  |  | 参与次数 | 运算能力标准分 |
|---|---|---|---|
| 参与次数 | Pearson 相关性 | 1 | 0.630** |
|  | 显著性（双侧） |  | 0.000 |
|  | N | 58 | 58 |
| 运算能力标准分 | Pearson 相关性 | 0.630** | 1 |
|  | 显著性（双侧） | 0.000 |  |
|  | N | 58 | 58 |

** 在 0.01 水平（双侧）上显著相关

从表 5-6A 和表 5-6B 中可知，实验组 A 和实验组 B 的学生参与讨论次数与其在调研中展现出的运算能力水平标准分在 0.01 水平（双侧）上显著相关，相关系数分别为 0.502 和 0.630，均达到了中等程度正相关的水平。由此可以得出，从学生层面来说，对自身运算能力的提升更需要主动有意识地参与，并持之以恒，才可能见到成效。

3. 影响学生参与讨论次数的相关性分析

在前述分析中已经可以看出，学生有意识地主动参与运算能力的讨论对于提升学生自身的运算能力水平有着正面的促进作用，但我们应该看到学生的这种主动有意识地参与与学生本身性别等固有因素没有特别明显的关系，对此笔者进行对应的相关性分析，如表 5-7 所示。

表 5-7 实验组学生参与讨论次数与性别的相关性分析

|  |  | 性别 |
|---|---|---|
| 实验组 A | Pearson 相关性 | 0.114 |
|  | 显著性（双侧） | 0.397 |
|  | N | 57 |
| 实验组 B | Pearson 相关性 | 0.016 |
|  | 显著性（双侧） | 0.906 |
|  | N | 58 |

如表 5-7 所示,实验组 A 和实验组 B 的学生参与讨论的次数与学生自身性别之间的 Pearson 相关系数分别为 0.114 和 0.016;显著性(双侧)水平分别为 0.397 和 0.906。这表明学生的性别差异与学生是否会更多地参与讨论之间并无显著相关性。这说明影响学生参与讨论的因素和性别没有关系,那么有什么因素会影响学生参与相关讨论的主动性呢?笔者对实验组学生的入围成绩标准分、参与讨论次数、调研成绩标准分进行了相关性分析,得到表 5-8。

表 5-8 实验组学生参与讨论次数、调研成绩与入围成绩的相关性分析

|  |  | 入围成绩标准分 |
|---|---|---|
| 参与讨论次数 | Pearson 相关性 | 0.433** |
|  | 显著性(双侧) | 0.000 |
|  | N | 115 |
| 调研成绩标准分 | Pearson 相关性 | 0.898** |
|  | 显著性(双侧) | 0.000 |
|  | N | 115 |

** 在 0.01 水平(双侧)上显著相关

从表 5-8 可以看出入围本次实践研究时学生的成绩标准分与其参与讨论次数之间的相关系数为 0.433,显著性水平为 0.000,在 0.01 水平(双侧)上显著相关,具有中等水平的相关性。而相关结果也表明这些学生入围时的成绩和参与期末调研时的成绩高度正相关,相关系数为 0.898。虽然不能对影响学生参与讨论次数的因素给出结论,但调研成绩对学生学习的主动性是有一定的正面推动作用的。

(三) 各评价维度分析
1. 期末调研运算能力内容维度结果分析

表 5-9 各组内容维度(知识板块)结果分析

| 分组 | 知识板块 | 总分 | 平均分 | 得分率 | 满分率 | 零分率 |
|---|---|---|---|---|---|---|
| 实验组 A | 方程与不等式 | 18 | 17.09 | 94.9% | 70.18% | 0% |
|  | 函数 | 13 | 8.42 | 64.8% | 5.26% | 0% |
| 实验组 B | 方程与不等式 | 18 | 13.86 | 77.0% | 32.76% | 3.45% |
|  | 函数 | 13 | 7.39 | 56.8% | 1.72% | 0% |
| 对照组 C | 方程与不等式 | 18 | 14.41 | 80.1% | 44.64% | 3.57% |
|  | 函数 | 13 | 7.89 | 60.7% | 1.79% | 0% |
| 整体 | 方程与不等式 | 18 | 15.12 | 84% | 49.12% | 2.34% |
|  | 函数 | 13 | 7.9 | 60.8% | 2.92% | 0% |

从上表 5-9 中得到的数据可以看出,两个知识板块的得分率最高达到了 84%,是方程与不等式部分,得分率较低的是函数部分的 60.8%,之所以存在较大的差距与两个知识板块间考察的数学运算素养水平有着直接的关系。其中方程与不等式所对应的题目为 13 题、17 题、18 题、22 题,共 18 分,其中只有 22 题第 2 问的 4 分对应的数学运算素养为水平二,其余 14 分均为运算素养水平一,反观函数部分涉及 1 题、19 题、26 题共 13 分,其中只有 1 题、19 题第 1 问、26 题第 1 问考察的数学运算素养水平为水平一共 6 分,19 题第 2 问和 26 题第 2 问对应数学运算素养水平二共 5 分,26 题第三问对应数学运算素养水平三占 2 分。由此可见九年级学生的数学运算素养水平大部分只能达到水平一的程度,经过一定的培养可以获得提高,但需要长期坚持培养才能见到比较明显的效果。

2. 期末调研运算能力素养水平维度结果分析

表 5-10 各组数学运算素养水平维度结果分析

| 分组 | 数学运算素养水平 | 总分 | 平均分 | 得分率 | 满分率 | 零分率 |
|---|---|---|---|---|---|---|
| 实验组 A | 水平一 | 20 | 19.3 | 96.6% | 68.42% | 0.00% |
|  | 水平二 | 9 | 5.9 | 65.7% | 7.02% | 1.75% |
|  | 水平三 | 2 | 0.3 | 14% | 10.53% | 82.46% |
| 实验组 B | 水平一 | 20 | 17.1 | 85.7% | 44.83% | 1.72% |
|  | 水平二 | 9 | 4.0 | 44.2% | 3.45% | 10.34% |
|  | 水平三 | 2 | 0.1 | 6.9% | 6.59% | 93.10% |
| 对照组 C | 水平一 | 20 | 17.4 | 87.2% | 50.00% | 0.00% |
|  | 水平二 | 9 | 4.7 | 52.4% | 3.57% | 10.71% |
|  | 水平三 | 2 | 0.1 | 7.1% | 7.14% | 92.86% |
| 整体 | 水平一 | 20 | 17.96 | 89.8% | 54.39% | 5.85% |
|  | 水平二 | 9 | 4.86 | 54.0% | 4.68% | 7.60% |
|  | 水平三 | 2 | 0.19 | 9.4% | 8.19% | 89.5% |

从表 5-10 可以看出，持续较长时间的专题讨论，对于提升学生数学运算素养水平一和水平二都取得了比较令人满意的效果，实验组 A 的学生在解答数学运算素养水平二的题目时满分率达到了 7.02%，明显高于实验组 B 和对照组 C，零分率只有 1.75%，远低于实验组 B 和对照组 C 的 10%，这表明通过对算理的分析和算法的对比，学生在处理数学运算素养水平二时有明显的提升

## 六、结论与展望
### （一） 研究结论

通过本次实践研究，可以得到以下结论：

1. 对数学运算能力的培养，需要尝试长时间的坚持，通过持续的对算理的分析和算法的对比，可以提升学生的数学运算素养水平。如果教师能够激发学生对数学运算的兴趣，让学生更主动地参与运算的相关学习研究，将能更好地提升学生的数学运算素养水平；
2. 通过相关性分析可以发现，学生参与实践讨论的次数越多与其数学运算素养水平的高低呈现出较为明显的正相关关系，相关系数在 0.502 到 0.630 之间。这表明学生越主动将精力投入到相关学习活动中，对于学生的提升越有益处；
3. 通过对运算能力内容维度的深入分析可以发现，九年级学生的数学运算素养水平仍处于较低状态，大部分学生经过学习可以达到水平一的程度，但在水平二和水平三层面上还需要进一步提高，这也应该是后续的高中学习中相关教师需要进一步予以关注的；
4. 学生数学运算能力水平与学生性别无显著相关性，但具有一定的持续性，呈现出强者恒强的趋势，所以数学运算作为六个核心素养之一，在解题中所处的位置决定了学生的数学学习中，具有较强的数学运算能力可以促进学生数学学业水平的提高。

### （二） 数学运算能力培养建议

根据本次实践研究的结果，笔者认为相关教师从以下方面着手，可以更高效地帮助学生提升其数学运算能力水平，培养其数学运算素养：

1. 激发学生学习运算的兴趣，调动学生的主观能动性，激发学生的参与意识，以典型问题适度激发学生的探讨热情，都有助于让学生主动参与数学运算的学习，提高学生学习运算的兴趣；
2. 数学运算能力的提升不能仅依靠练习题目，还需要在选题上多下功夫，让题目更具有典型性，可以帮助学生更好地明晰算理；
3. 在明晰算理的同时要讲究方法着力培养学生良好的数学思维品质，通过分析题目选择更优的运算思路，设计合理的运算程序，提升学生的算法选择能力。

### 七、研究反思与展望

本次实践研究中还存在很多问题，可以在之后的相关研究中进行改进。

1. 实验组和对照组设置有待改进，目前的设置中后半个学期只剩余一个实验组A，应该设立后半学期参与的另一个实验组，以完备整个实践的设计；
2. 本次实践研究中对学生数学运算能力水平的评价数据收集不足，仅收集对比了八年级第二学期和九年级第一学期各两次大考的调研结果作为评价数据，缺乏实践过程中的其他数据支持，无法观察到学生数学运算能力水平提升的过程。需要在未来的研究中加以改进，收集更多的过程性数据，比如在讨论课中设置相应的反馈环节，及时收集过程性数据；
3. 使用区调研成绩作为评价数据的来源虽然可靠性得到了保证，但限于教学进度，这些调研无法全面反应学生的数学运算素养水平，毕竟九年级第一学期期末的调研重点在于方程与不等式和函数部分，对于数与式部分的关注略有不足。而且此次区调研更多的关注了数学运算素养水平一和水平二的考察，对水平三的考察做的并不到位，同时调研中其他方向的考题也分散了学生的精力，学生无法专注于完成运算部分的调研。鉴于此，仍有必要在未来编制专门用于测量学生数学运算素养水平的调研试卷；
4. 本次研究更多的关注了技术层面去提升学生的数学运算能力水平，未能有效涉及学生的情感态度价值观维度，需要在今后的其他研究中进行完善。

### 八、参考文献

[1] 教育部. 义务教育数学课程标准（2011年版）[M].北京：北京师范大学出版社，2011
[2] 教育部. 普通高中数学课程标准（2017年版2020年修订）[M].北京：人民教育出版社，2020
[3] 黄友初. 数学素养的内涵、测评与发展研究[M].北京：科学出版社，2016
[4] 北京市海淀区教师进修学校. 海淀区义务教育学业标准与教学指导[M].北京：北京师范大学出版社，2018
[5] 张庆辰. 新课标下初中生数学运算能力的培养策略研究[D].延边大学，2019.04

**九、附录**

**附录 1：运算素养培养专题讨论示例（第 6 课时）**

**一、教学目标：**
1、能用不同的方式解决有关不等式问题中的比大小问题；
2、能有意识地将有关不等式问题中的比大小问题转化为函数问题；
3、能借助函数观点解决有关不等式中的比大小问题；
4、在解决比大小问题的过程中，能根据实际情况选择合适的思路解决问题；
5、感悟数学知识之间的联系，认识函数的重要性，进一步体会化归与转化思想、数形结合思想、分类讨论等数学思想，从中提升数学抽象、逻辑推理和数学运算素养。

**二、教学重、难点：**
重点：不同比较代数式大小的思路的比较
难点：能将有关不等式问题中的比大小问题转化为函数问题.

**三、典例分析：**
例 1：若 $a > b$，比较 $-3.5b + 1$ 与 $-3.5a + 1$ 的大小关系；
解：① 借助不等式的性质比大小

$\because a > b$

$\therefore -3.5a < -3.5b$ （不等式性质 3）

$\therefore -3.5a + 1 < -3.5b + 1$ （不等式性质 1）

② 借助做差和 0 比大小

$-3.5a + 1 - (-3.5b + 1)$

$= -3.5a + 1 + 3.5b - 1$

$= -3.5a + 3.5b$

$= 3.5(b - a)$

$\because a > b$

$\therefore b - a < 0$

$\therefore 原式 = 3.5(b - a) < 0$ （乘法法则）

故 $-3.5a + 1 < -3.5b + 1$

③ 借助函数观点比大小

令 $y = -3.5x + 1$

则当 $x = a$ 时，有 $y_a = -3.5a + 1$；

当 $x = b$ 时，有 $y_b = -3.5b + 1$

$\because 一次函数 y = -3.5x + 1$ 中，$y$ 随 $x$ 增大而减小

又 $\because a > b$

$\therefore y_a < y_b$

即 $-3.5a + 1 < -3.5b + 1$

【设计意图】展现不同方法比大小时的特点进行对比，锻炼学生基本的运算素养，特别是当采用"不等式的性质变形"和"做差与 0 比大小"两种方法时，要让学生明确其中每一步所

依据的算理，将对学生的数学运算能力与逻辑推理能力的培养有所帮助，而对比借助函数图象比大小的思路，其中对数形结合的运用，也会为学生打开思路，有利于后续问题的探究。

例2：（2019年北京中考）用三个不等式 $a>b$，$ab>0$，$\frac{1}{a}<\frac{1}{b}$ 中的两个不等式作为题设，余下的一个不等式作为结论组成一个命题，组成真命题的个数为（　　）

A．0　　　　　　B．1　　　　　　C．2　　　　　　D．3

解：①若 $a>b$，$ab>0$，则 $\frac{1}{a}<\frac{1}{b}$；真命题：

理由：$\because a>b$，$ab>0$，

$\therefore \frac{a}{ab} > \frac{b}{ab}$ （不等式性质2）

$\therefore \frac{1}{a} < \frac{1}{b}$；

②若 $ab>0$，$\frac{1}{a} < \frac{1}{b}$，则 $a>b$，真命题；

理由：$\because ab>0$，

$\therefore \frac{1}{a} \times ab < \frac{1}{b} \times ab$，（不等式性质2）

$\therefore a>b$．

③若 $a>b$，$\frac{1}{a} < \frac{1}{b}$，则 $ab>0$，真命题；

理由：$\because \frac{1}{a} < \frac{1}{b}$，

$\therefore \frac{1}{a} - \frac{1}{b} < 0$，即 $\frac{b-a}{ab} < 0$，

$\because a>b$，

$\therefore b-a < 0$，

$\therefore ab>0$ （乘法法则）

$\therefore$ 组成真命题的个数为3个；

解法2：令 $y = \frac{1}{x}$

则当 $x=a$ 时，有 $y_a = \frac{1}{a}$；

当 $x=b$ 时，有 $y_b = \frac{1}{b}$

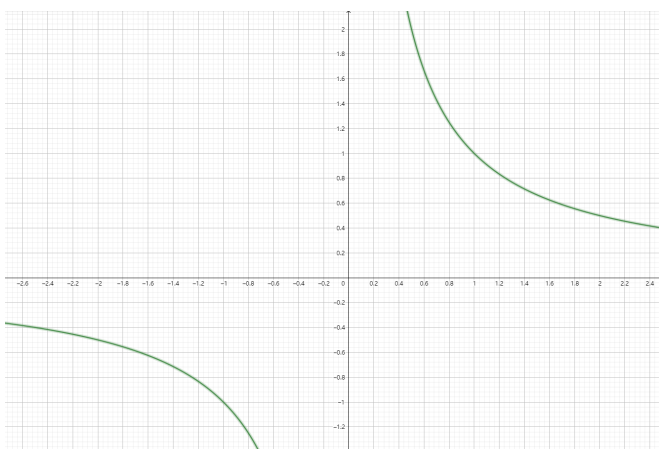

由图象容易得出，

① 若 $a>b$，$ab>0$，则 $\frac{1}{a} < \frac{1}{b}$，是真命题；

② 若 $ab>0$，$\frac{1}{a} < \frac{1}{b}$，则 $a>b$，是真命题；

③ 若 $a>b$，$\frac{1}{a} < \frac{1}{b}$，则 $ab>0$，是真命题；

【设计意图】学生通过利用不等式性质和作差法比大小进一步熟悉的锻炼数学运算能力，之后再通过函数观点来进行判断，通过函数图象，借助数形结合的思想，通过简单的逻辑推理就可以更容易地判断命题的真假，由此学生可以发现不同的角度在解题时运算量的大小和难度是不同的，需要从不同的角度出发去尝试解决问题，并在尝试的过程中发现相对更优的运算思路。

例3：如果 $a$，$b \in R$，且 $a > b$，试比较：$a^2$ 与 $b^2$ 的大小；

解法1：作差比大小

$a^2 - b^2 = (a+b)(a-b)$

$\because a > b$

$\therefore a - b > 0$

$\therefore a^2 - b^2$ 的正负性与 $a + b$ 的正负性保持一致

当 $a > b > 0$ 时，$a + b > 0$，故 $a^2 - b^2 > 0$，即 $a^2 > b^2$；

当 $a > 0 > b$ 时，根据加法法则，$a + b$ 的正负性由 $a$ 与 $|b| = -b$ 的大小关系决定

　　当 $a > -b$ 时，$a + b > 0$，故 $a^2 - b^2 > 0$，即 $a^2 > b^2$；

　　当 $a = -b$ 时，$a + b = 0$，故 $a^2 - b^2 = 0$，即 $a^2 = b^2$；

　　当 $0 < a < -b$ 时，$a + b < 0$，故 $a^2 - b^2 < 0$，即 $a^2 < b^2$；

当 $0 \geq a > b$ 时，$a + b < 0$，故 $a^2 - b^2 < 0$，即 $a^2 < b^2$；

综上所述，

当 $a > b > 0$ 时，$a + b > 0$，故 $a^2 - b^2 > 0$，即 $a^2 > b^2$；

当 $b < 0$ 且 $a > -b$ 时，$a + b > 0$，故 $a^2 - b^2 > 0$，即 $a^2 > b^2$；

当 $b < 0$ 且 $a = -b$ 时，$a + b = 0$，故 $a^2 - b^2 = 0$，即 $a^2 = b^2$；

当 $b < 0$ 且 $b < a < -b$ 时，$a + b < 0$，故 $a^2 - b^2 < 0$，即 $a^2 < b^2$。

解法2：令 $y = x^2$，

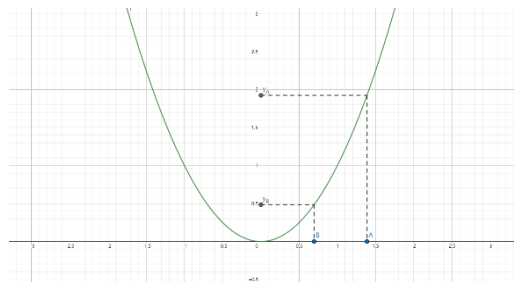

　　则当 $x = a$ 时，有 $y_a = a^2$；

　　当 $x = b$ 时，有 $y_b = b^2$

　　由二次函数图象可知，对称轴为 $x = 0$

　　$\because$ 二次项系数大于 0，

　　$\therefore$ 在对称轴两侧，$y$ 随 $|x|$ 的增大而增大

　　$\because a > b$

　　$\therefore$ 当 $a > |b|$ 时，$a + b > 0$，故 $a^2 - b^2 > 0$，即 $a^2 > b^2$；

　　当 $b < 0$ 且 $a = -b$ 时，$a + b = 0$，故 $a^2 - b^2 = 0$，即 $a^2 = b^2$；

　　当 $b < 0$ 且 $a < |b|$ 时，$a + b < 0$，故 $a^2 - b^2 < 0$，即 $a^2 < b^2$。

【设计意图】学生首先会利用作差法比大小，这个过程中锻炼学生的数学运算能力，有助于提升其对分类讨论思想的认识和数学运算素养. 基于作差法比较繁琐，教师引导学生通过观察式子结构特征，思考一个更为直观的做法，让学生学会将有关不等式问题中的比大小问题分别转化为反比例函数、二次函数等函数的单调性问题，体会化归与转化思想、数形结合的

思想，同时提升了逻辑推理、数学抽象素养.

## 四、反思提升（本节课你的收获与不足）

1. 思考：如何比较两数（代数式）的大小？（本节课涉及的方法）
（1）利用不等式的性质进行变换；
（2）作差：作差后通过分解因式、配方等手段判断差的符号得出结果；
（3）利用函数的单调性及图象法.
其中比较法（作差、作商）是最基本的方法. 本节课根据题目式子的结构特征我们可以抽象出对应函数，利用函数的单调性、对称性（二次函数）等性质比较大小，简单直观，不易错！

## 总结：

  本节课我们学习了从不同视角去解不等式比大小问题的基本思路和方法，通过对比我们可以发现在不同的条件下，选择适当地思路去进行运算，其复杂程度和推理运算的难度有着明显的差异，在这个过程中可以深刻体会化归与转化思想、数形结合思想、分类讨论思想等一系列数学思想方法。

  特别是在一定条件下，借助联想并构造合适的函数，从函数观点去看不等式中的比大小问题可以明显起到简化运算的作用，进而通过函数的性质可以相对更简单的解决有关问题。
【设计意图】巩固本节课内容—利用不同的运算思路去解决代数式比大小的问题，特别是将有关不等式问题中的比大小问题转化为函数问题在特定条件下可以有效简化运算，希望学生逐渐养成反思的好习惯，并学会在反思中成长进步。

## 五、检测反馈：

1、已知点（$\sqrt{2}+1$，$y_1$）和（$\sqrt{3}-1$，$y_2$）在二次函数$y = x^2 - 2x + a$的图象上，试比较$y_1$和$y_2$的大小.
2、如果$a, b \in R$，且$a > b,$试用不同的方法比较：
【设计意图】了解学生课上的学习成果.

**附录 2：2021-2022 学年度学区第一学期期末学业水平调研九年级数学试卷**

仅包括表 5-1 中列出本次实践研究使用的题目

1. 在平面直角坐标系 $xOy$ 中，下列函数的图象经过点 $(0,0)$ 的是

   (A) $y = x+1$  (B) $y = x^2$  (C) $y = (x-4)^2$  (D) $y = \dfrac{1}{x}$

13. 若关于 $x$ 的方程 $x^2 - 2x + k = 0$ 有两个不相等的实数根，则 $k$ 的取值范围是________．

17. 解方程：$x^2 - 6x + 8 = 0$．

18. 已知 $a$ 是方程 $2x^2 - 7x - 1 = 0$ 的一个根，求代数式 $a(2a-7)+5$ 的值．

19. 在平面直角坐标系 $xOy$ 中，抛物线 $y = a(x-3)^2 - 1$ 经过点 $(2,1)$．
    （1）求该抛物线的表达式；
    （2）将该抛物线向上平移________个单位后，所得抛物线与 $x$ 轴只有一个公共点．

22. 已知关于 $x$ 的一元二次方程 $x^2 + (2-m)x + 1 - m = 0$．

    （1）求证：方程总有两个实数根；
    （2）若 $m < 0$，且此方程的两个实数根的差为 $3$，求 $m$ 的值．

26. 在平面直角坐标系 $xOy$ 中，点 $(4,3)$ 在抛物线 $y = ax^2 + bx + 3(a>0)$ 上．

    （1）求该抛物线的对称轴；
    （2）已知 $m>0$，当 $2-m \leq x \leq 2+2m$ 时，$y$ 的取值范围是 $-1 \leq y \leq 3$，求 $a$，$m$ 的值；
    （3）在（2）的条件下，是否存在实数 $n$，当 $n-2 < x < n$ 时，$y$ 的取值范围是 $3n-3 < y < 3n+5$，若存在，直接写出 $n$ 的值；若不存在，请说明理由．